# Modeling a demographic problem using the Leslie matrix

O.A. Malafeyev
Saint-Petersburg State University
o.malafeev@spbu.ru

T.R. Nabiev
Saint-Petersburg State University
obiwankenobi731@gmail.com

N.D. Redinskikh
Saint-Petersburg State University
redinskich@yandex.ru

**Key words:** geopolitics, two person zero-sum game, preferences, Leslie matrix, Nash equilibrium, demographics, immigration, game theory, age-classes, survival table, normal-form game.

## Abstract

The application of Leslie matrices in demographic research is considered in this paper. The Leslie matrix is first proposed in the 1940s and gained popularity in the mid-1960s, becoming fundamental tool for predicting population dynamics. The Leslie matrix allows to categorize individuals based on various attributes and calculate the expected population sizes for various demographic categories in subsequent time intervals. The universality of the Leslie matrix extends to diverse life cycles in plants and animals, making it ubiquitous tool in non-human species. In the paper is presented detailed application of Leslie matrices to the problem of the two countries, demonstrating their practical value in solving real demographic problems. In conclusion, the Leslie matrix remains a cornerstone of demographic analysis, reflecting the complexity of population dynamics and providing a robust framework for understanding the intricate interplay of factors shaping human society. Its enduring relevance and adaptability make it an essential component in the toolkit of demographers and ecologists.

# Introduction

Demography, the statistical study of the population, is an area that provides a critical understanding of the structure and dynamics of human society. Encompassing variables like birth rate, death rate, age distribution, and migration patterns, demography relies significantly on analytical tools such as the Leslie matrix. This mathematical construct proves instrumental in modeling changes within population demographics over time, providing a comprehensive framework for demographic studies that transcend simplistic life and death considerations, incorporating a spectrum of individual attributes.

The Leslie Matrix is named after ecologist Patrick H. Leslie, who developed it in the middle of the 20th century. Leslie Matrix represents a demographic projection matrix that categorizes individuals based on age, gender, employment status, maturity, and other pertinent factors. Such forecast takes into account the probability of people moving from one category to another, for example, from the unemployed to the employed. By incorporating these probabilities, the Leslie matrix facilitates the computation of anticipated population sizes for distinct demographic categories at subsequent temporal intervals.

In this paper, the structures of population matrix models classified by age and closely related to Leslie matrices are considered. This also applies to population groups classified by other factors, which emphasizes the universality of matrix methods in demographic analysis. The Leslie matrix is especially important to biologists because of the diversity of plant and animal life cycles.

The application of Leslie matrices in demographic research offers a reliable basis for predicting population dynamics. In this paper it is studied the application of Leslie matrices in demographic tasks. The mathematical basis of the Leslie matrix, its practical application and the impact of demographic trends on policy and planning are considered. The methods and approaches presented in [70-120] are also used in this paper.

# Leslie matrices

The Leslie matrix, a tool employed in demographic investigations, was initially proposed in the 1940s by Bernardelli, Lewis, and notably Leslie, after whom it is named. However, it garnered minimal attention until the mid-1960s when ecologists and demographers rediscovered its utility, marking the onset of its prominence in demographic research. The Leslie matrix serves to predict population dynamics, considering variables such as age, gender, employment status, and maturity. Its application enhances the comprehension of demographic changes, surpassing the simplistic life-and-death dichotomy.

Within each temporal interval, individuals of each type undergo the probabilities of death, transition to alternative categories (e.g., employment for the unemployed), or reproduction. By incorporating these probabilities, the anticipated population of each type in the subsequent interval is calculated, constituting a population forecast. This forecast is conveniently represented by a population forecast matrix, facilitating the computation of various population dynamics characteristics.

Leslie was an ecologist with a particular interest in small mammal populations (Crowcroft 1991). In addition to his work on matrix population models, he calculated for the first time in the survival table the internal growth rate for all species except humans, and made significant contributions to stochastic models and the assessment of recapture. Matrix models were mostly neglected until the mid-1960s, when they are rediscovered by both ecologists (Lefkovich, 1965) and demographers (Lopez, 1961, Keifitz, 1964, Rogers, 1966).

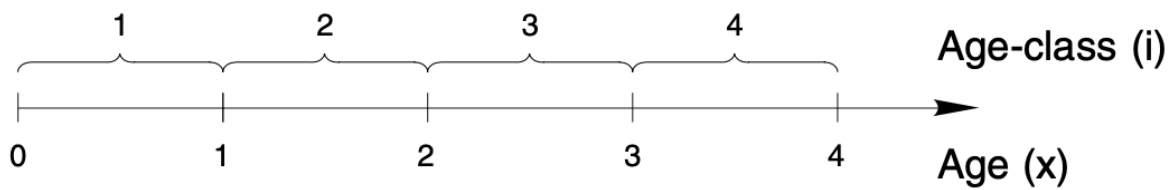

Figure 1. The relationship between continuous age variable and discrete age classes.

Commence with a model delineating the population by age. Discretize the continuous age variable, commencing from 0, into distinct age classes originating from 1. The schematic representation is presented in Figure 1. Age class i corresponds to the age range i − 1 ≤ x ≤ i. According to this agreement, the first age class is numbered as 1. It is noteworthy that certain authors may opt to designate the first age class as 0.

The objective is to forecast the population from time t to time t+1. Let us assume that the unit of time equals the width of the age class. This unit is referred to the projection interval, and its selection is one of the initial steps in constructing a matrix model. It is unsurprising that a model forecasting year by year will differ from one forecasting month by month or decade by decade.

Let us assume that the forecasting interval is one year, and individuals are categorized into three age groups (0–1, 1–2, and 2–3 years). The population's state is described by the vector n(t), where its elements $n_i(t)$ indicate the quantity of individuals in each age class. Individuals in age classes 2 and 3 at time t + 1 are individuals who have survived the preceding age classes at time t.

$$n_2(t+1) = P_1 n_1(t)$$
$$n_3(t+1) = P_2 n_2(t),$$

Here, $P_i$ represents the probability that an individual of age class i will survive for one year. New members of age class 1 cannot originate from any other age class; they must have resulted from reproduction. Thus, we express

$$n_1(t+1) = F_1 n_1(t) + F_2 n_2(t) + F_3 n_3(t),$$

where $F_i$ denotes the birth rate per capita * of the population in age class *i*, i.e. the number of individuals in age class 1 at time t+1 per person in age class i at time t.

These equations are conveniently expressed in matrix form as

$$\begin{pmatrix} n_1 \\ n_2 \\ n_3 \end{pmatrix}(t+1) = \begin{pmatrix} F_1 & F_2 & F_3 \\ P_1 & 0 & 0 \\ 0 & P_2 & 0 \end{pmatrix} \begin{pmatrix} n_1 \\ n_2 \\ n_3 \end{pmatrix}(t)$$

Or, more succinctly,

$$\mathbf{n}(t+1) = \mathbf{A}\mathbf{n}(t),$$

where **n** is the population size vector, and **A** is the population forecast matrix. This specific version with age classification is often referred to as the Leslie matrix. It is non-negative (as negative elements imply a negative number of individuals), with positive entries only in the first row and sub-diagonal (representing survival probabilities).

Population matrix models can be classified based on the nature of the projection matrix **A**. In the simplest case, the matrix remains constant.

## The Two Countries Problem

Let's consider the case of two countries, A and B. Both nations have populations that is categorized into three groups: children, adults, and elderly. For each country, the initial population, expressed in thousands, is represented by vectors:

$$A = \begin{pmatrix} 30 \\ 35 \\ 25 \end{pmatrix}, B = \begin{pmatrix} 40 \\ 30 \\ 30 \end{pmatrix},$$

For each country, the population for the next unit of time is determined by groups as follows:

$$\begin{pmatrix} n_1 \\ n_2 \\ n_3 \end{pmatrix}(t+1) = \begin{pmatrix} F_1 & F_2 & F_3 \\ P_1 & 0 & 0 \\ 0 & P_2 & 0 \end{pmatrix} \begin{pmatrix} n_1 \\ n_2 \\ n_3 \end{pmatrix}(t) + \begin{pmatrix} i_1 \\ i_2 \\ i_3 \end{pmatrix}$$

Where the last vector represents an additional vector of immigrants who have migrated to this country according to age groups.

For each country, there are initial Leslie matrices and immigration vectors.

$$LA_0 = \begin{pmatrix} 0 & 2 & 1 \\ 0.2 & 0 & 0 \\ 0 & 0.4 & 0 \end{pmatrix} \quad IA_0 = \begin{pmatrix} 5 \\ 10 \\ 10 \end{pmatrix}$$

$$LB_0 = \begin{pmatrix} 0 & 5 & 2 \\ 0.2 & 0 & 0 \\ 0 & 0.4 & 0 \end{pmatrix} \quad IB_0 = \begin{pmatrix} 15 \\ 20 \\ 20 \end{pmatrix}$$

Both countries aim to improve their demographics, and each state has two strategies: to develop healthcare within the country, reducing mortality in each age group, or to enhance the quality of life and provide incentives for immigrants, increasing the number of immigrants to the country. Consequently, if there is an investment in survival (S), the corresponding Leslie matrices will appear as follows:

$$LA_S = \begin{pmatrix} 0 & 3 & 1 \\ 0.4 & 0 & 0 \\ 0 & 0.6 & 0 \end{pmatrix} \quad LB_S = \begin{pmatrix} 0 & 6 & 2 \\ 0.6 & 0 & 0 \\ 0 & 0.8 & 0 \end{pmatrix}$$

In this case, the investment outcome is independent of the tactic chosen by the competing country, unlike the strategy where one or both countries opt for immigration investments. If country A invests in immigration, a portion of the population immigrating to this country emigrates from country B, or chooses country A for emigration instead of country B. In this scenario, it can be represented as an increase in the number of immigrants to country A and a decrease in the number of emigrants to country B, as illustrated in the following:

$$IA_i = \begin{pmatrix} 35 \\ 40 \\ 40 \end{pmatrix} \quad IB_{-i} = \begin{pmatrix} 10 \\ 15 \\ 15 \end{pmatrix}$$

Similarly, if country B chooses to invest in immigration while country A invests in survival, this can be expressed as follows:

$$IB_i = \begin{pmatrix} 45 \\ 50 \\ 50 \end{pmatrix} \quad IA_{-i} = \begin{pmatrix} -5 \\ 0 \\ 0 \end{pmatrix}$$

If both countries choose to invest in immigration, the vector $I$ for each country is remained unchanged.

If both countries, A and B, want to maximize the population in each age group, an important question arises: who decides what to do, country A or country B? There are two scenarios: country A decides first, or country B decides first.

Let us examine each case sequentially. Let us assume that country A decides first, such situation is represented on Figure 2.

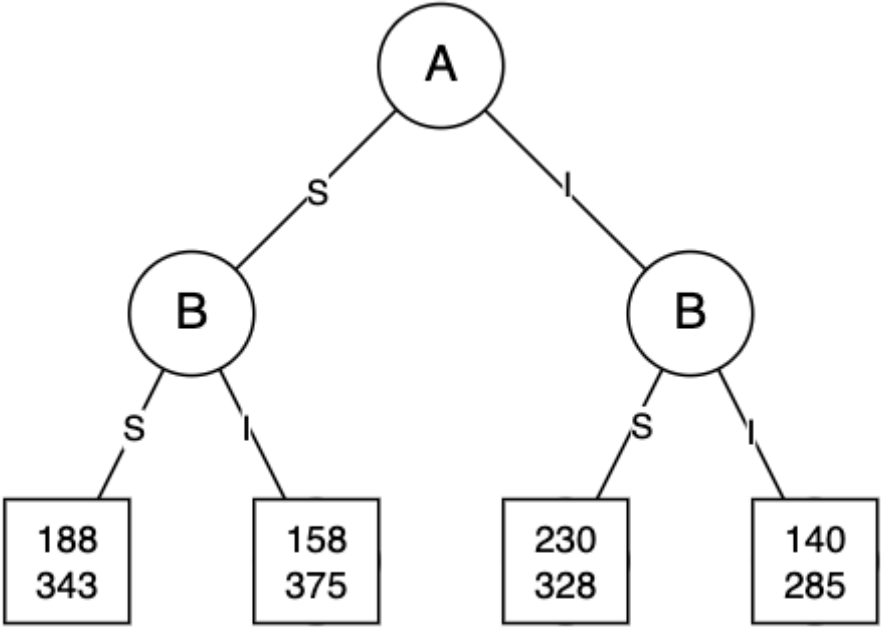

Figure 2. The game diagram where country A makes the first move.

In the game tree diagram, the top part of the tree contains the root node (a circle labeled A) with two branches, labeled S (survival) and I (immigration). This indicates that country A can choose to go left (S) or right (I). This leads us to two nodes labeled B, in each of which country B can invest in survival (S) or immigration (I).

While country A has only two strategies, country B has four:
1. Invest in survival regardless of country A's move (SS).
2. Invest in immigration regardless of country A's move (II).
3. Mimic what country A does (SI).
4. Do the opposite of what country A does (IS).

The first letter in parentheses indicates country B's move if country A invests in survival, and the second letter indicates country B's move if country A invests in immigration.

The move taken by a player in a node is an action, and a series of actions that fully define a player's behavior is a pure strategy. Thus, country A has two strategies, each representing a single action, while country B has four strategies, each representing two actions; one will be used when country B goes left, and the other when country B goes right.

At the bottom of the game tree, there are four nodes known as leaf or terminal nodes. Each terminal node contains a payoff for both players: first for country A (player 1), and then for country B (player 2), if they choose strategies that lead them to that specific leaf. How should country A decide what to do?

Clearly, country A needs to figure out how country B reacts to each of the two choices by country A: S and I. If country A chooses S, country B chooses I because it results in a greater population increase of 375. Thus, country A gains 158 by moving right. If country A chooses I, country B chooses S, as selecting S yields 328 instead of 285. Therefore, country A gains 158 for choosing S and 230 for choosing I. So country A should choose I.

A Nash equilibrium in a two-player game is a pair of strategies, each of which is the best response to the other; in other words, each provides the player using it with the maximum possible payoff, considering the strategy of the other player. Let us represent this game in normal form. Both representations are commonly used and can be switched depending on convenience. The normal form corresponding to Figure 2 is presented in Table 1.

|  |  | country B |  |  |  |
|---|---|---|---|---|---|
|  |  | SS | SI | IS | II |
| country A | S | 188,343 | 188,343 | 158,375 | 158,375 |
|  | I | 230,328 | 140,285 | 230,328 | 140,285 |

Table 1. Normal-form game

In this example, player 1's (Country A) strategies are arranged in rows, and player 2's (Country B) strategies are in columns. Each entry in the resulting matrix represents the payoffs for both players if they choose the corresponding strategies. To find Nash equilibria from the normal form of the game, let us attempt to select a row and column such that the payoff at their intersection is the maximum possible for player 1 in the column and the maximum possible for player 2 along the row (there may be multiple such pairs). In this case, there are two Nash equilibria: strategies (I, SS) and (I, IS).

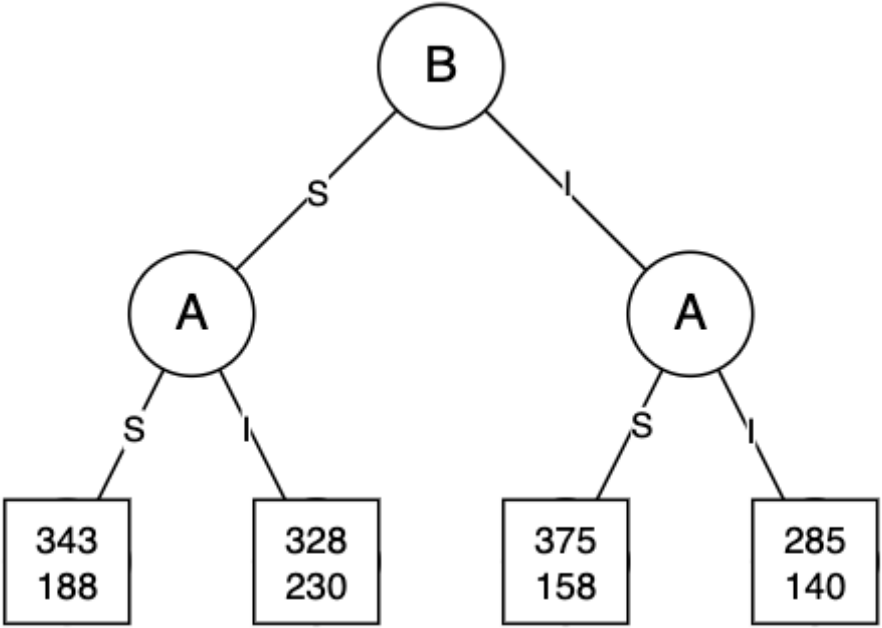

Figure 3. Game Tree where Country B Makes the First Move

But what if Country B has to choose first? This is a game with an extensive form (Figure 3). Now, we label Country B as player 1 and Country A as player 2. Now, Country A has four strategies (the strategies that belonged to Country B in the previous version of the game), and Country B has only two (those that belonged to Country A earlier).

Country B observes that the best response of Country A to S is I, and the best response of Country A to I is S. Since Country B gets 328 in the first case and only 375 in the second, Country B chooses I. Then the optimal choice for Country A will be IS, and the payoffs will be (375, 158).

The normal form of the case where Country B moves first is shown in Table 2. Once again, we find two Nash equilibria (I; SS) and (I; IS).

|  | Country A | | | |
|---|---|---|---|---|
|  | SS | SI | IS | II |

| Country B | S | 343,188 | 343,188 | 328,23 | 328,23 |
|---|---|---|---|---|---|
| | I | 375,158 | 285,140 | 375,158 | 285,140 |

Table 2. Normal-form game

## Conclusion

In conclusion, it is noteworthy that the Leslie matrix has proven to be an indispensable tool in the field of demography, offering sophisticated means of understanding and predicting population dynamics. Utilizing this matrix, demographers and ecologists can model populations with a level of detail that takes into account various individual characteristics, such as age, gender, employment status, and maturity.

The versatility of the Leslie matrix allows for the consideration of complex life cycle phenomena, particularly relevant in the study of human demography. For human populations, the matrix adapts to include changes in employment status and other socio-economic factors, providing a more detailed representation of population changes than traditional mortality tables.

As revealed in this example, the Leslie matrix is not merely a theoretical construct but a practical tool applied to solve real demographic challenges. In this case, it facilitated the development of policies and planning strategies, enabling accurate forecasting of demographic trends. The matrix's ability to classify individuals based on various characteristics and model probabilities of transitioning between these states has made it a cornerstone of demographic analysis.

In summary, the Leslie matrix remains a cornerstone of demographic analysis, reflecting the complexity of population dynamics and providing the foundation for understanding the intricate interplay of factors shaping human society. Its enduring relevance attests to its reliability and adaptability, making it a crucial component in the toolkit of both demographers and ecologists.

preobrazovanij XXI veka» // V sbornike: Filosofiya poznaniya i tvorchestvo zhizni. Sbornik statej – 2014 - S. 279-292.